\documentclass{ifacconf}
\usepackage{graphicx}      
\usepackage{natbib}        
\usepackage{amsmath,amssymb,amsfonts}

\usepackage{algorithm}
\usepackage{algorithmicx}
\usepackage{algpseudocode}

\usepackage{textcomp}
\usepackage{xcolor}
\usepackage{rotating}
\usepackage{mathrsfs}
\usepackage{tikz}
\usetikzlibrary{arrows} 
\newcommand{\bmat}{\left[ \begin{matrix}}
	\newcommand{\emat}{\end{matrix} \right]}
\newcommand{\innerprod}[2]{\langle{#1},\,{#2}\rangle}\label{key}
\DeclareMathOperator{\E}{{\mathbb E}}
\newcommand{\Rbb}{\mathbb R}

\newcommand{\Cbb}{\mathbb C}

\newcommand{\Zbb}{\mathbb Z}

\newcommand{\Tbb}{\mathbb T}

\newcommand{\qb}{\mathbf q}
\newcommand{\tb}{\mathbf t}
\newcommand{\kb}{\mathbf k}
\newcommand{\lb}{\boldsymbol{\ell}}

\newcommand{\zerob}{\mathbf 0}

\newcommand{\Nb}{\mathbf N}

\newcommand{\Qb}{\mathbf Q}
\newcommand{\Rb}{\mathbf R}

\newcommand{\thetab}{\boldsymbol{\theta}}

\newcommand{\Sigmab}{\boldsymbol{\Sigma}}

\newcommand{\vtr}{\mathrm{vec}}

\newcommand{\Lscr}{\mathscr{L}}

\newcommand{\Hcal}{\mathcal{H}}
\newcommand{\Gcal}{\mathcal{G}}
\newcommand{\Wcal}{\mathcal{W}}
\renewcommand{\d}{\mathrm{d}}

\begin{document}
\begin{frontmatter}

\title{\LARGE \bf
	Revisiting a Fast Newton Solver for a 2-D Spectral Estimation Problem: Computations with the Full Hessian} 

\thanks[footnoteinfo]{This work was supported in part by the National Natural Science Foundation of China (No.~62103453) and the Innovation Training Program for College Students of Sun Yat-sen University (No.~20252496).}

\author[First]{Ji Cheng and Bin Zhu} 

\address[First]{School of Intelligent Systems Engineering, Sun Yat-sen University, Gongchang Road 66, 518107 Shenzhen, China (e-mails: chengj227@mail2.sysu.edu.cn, zhub26@mail.sysu.edu.cn)}

\begin{abstract}                

Spectral estimation plays a fundamental role in frequency-domain identification and related signal processing problems.
This paper revisits a 2-D spectral estimation problem formulated in terms of convex optimization.
More precisely, we work with the dual optimization problem and show that the full Hessian of the dual function admits a Toeplitz-block Toeplitz structure which is consistent with our finding in a previous work.
This particular structure of the Hessian enables a fast inversion algorithm in the solution of the dual optimization problem via Newton's method
whose superior speed of convergence is illustrated via simulations.

\end{abstract}


\end{frontmatter}

\section{Introduction}

Spectral estimation is about inferring the statistical power distribution of a random process/field in the frequency domain from a finite number of samples. It has vast applications in e.g., system identification \citep{ringh2015multidimensional,ringh2018multidimensional,Zhu-Zorzi-2022-cep-est}, radar signal processing \citep{engels2017advances,ZFKZ2019fusion,ZFKZ2019M2}, texture image generation \citep{lindquist2016modeling}, and geotechnical modeling \citep{Liu-Zhu-2021-rand-field}.

Among different methods for spectral estimation, there is a series of works that center around the idea of \emph{rational covariance extension}, see e.g., \cite{BGL-98,BGL-THREE-00,SIGEST-01,Georgiou-L-03,FPR-08,FMP-12,Z-14rat,Z-14,ZFKZ2020M2-SIAM,Zhu-Zorzi-2023-cepstral,Zhu-Zorzi-2023-entro-inc}. In those works, spectral estimation is formulated as equality constrained optimization problems in which the constraints are \emph{trigonometric moment equations} for the spectral density, and the objective functional is usually some kind of entropy or a related concept of divergence. These optimization problems are successful because they can be made well-posed and allow for solution tuning by smoothly changing the parameters.

In a recent work \cite{Liu-Zhu-2021}, a spectral estimation problem for two-dimensional (2-D) random fields was considered using the \emph{Itakura-Saito} (IS) pseudo-distance \citep{FMP-12} as the objective functional in the optimization. The problem was solved via duality, and the Hessian of the dual function was shown to have a \emph{Toeplitz-block Toeplitz} (TBT) structure \citep{wax1983efficient} which makes possible a fast Newton algorithm. However, only a principal submatrix of the full Hessian was computed  in that paper, and thus strictly speaking, the algorithm should be interpreted as a \emph{quasi-Newton} method \citep{kreutz2009complex} which in general, converges slower than the true Newton's method. This last observation motivates the investigation in this paper. 
Indeed, we show that the TBT structure of the full Hessian is preserved if we collect the dual variables in a suitable order. Hence the fast inversion algorithm in \cite{wax1983efficient} can be used provided that proper modifications are made for Hermitian matrices. In addition, we give numerical examples in which Newton's method, enabled by the fast inversion algorithm for the full Hessian, converges much faster than the quasi-Newton method in \cite{Liu-Zhu-2021}.

The remaining part of this paper is organized as follows.
Sec.~\ref{sec:rev_prob} gives a brief review of a 2-D spectral estimation problem formulated in terms of optimization.
Sec.~\ref{sec:struct_Hessian} characterizes the algebraic structure of the full Hessian of the dual objective function.
Sec.~\ref{sec:fast_inv} describes a fast inversion algorithm for the structured Hessian, and Sec.~\ref{sec:TBT_linear_eq} treats the computation of the Newton direction, i.e., solving linear equations where the Hessian is the coefficient matrix.
Sec.~\ref{sec:sims} provides simulation results.
Sec.~\ref{sec:conclus} concludes the paper.

\section{Review of a 2-D spectral estimation problem}\label{sec:rev_prob}


We focus on the spectral estimation problem for complex-valued zero-mean wide-sense stationary $2$-D random field $\{y_{\tb} : \tb=(t_1, t_2) \in \Zbb^2\}$. The assumption of stationarity implies that the (auto-)covariance of the random field $\sigma_{\kb}:=\E y_{\tb+\kb}\overline{y_{\tb}}$ is a function of the difference $\mathbf{k}\in\Zbb^2$ between the vectorial indices only. Here $\E$ denotes the mathematical expectation and $\overline{z}$ signifies the complex conjugate of $z\in\Cbb$.
The (power) spectral density of the random field is obtained by taking the Fourier transform of the covariance function $\sigma_{\kb}$, i.e.,
\begin{equation}\label{psd_def}
	\Phi(\thetab) := \sum_{\mathbf{k}\in \Zbb^2}\sigma_\mathbf{k}e^{-\imath\langle\mathbf{k},\boldsymbol{\theta}\rangle},
\end{equation}
where $\thetab=(\theta_1,\theta_2)$ such that each $\theta_j,\ j=1, 2$ is an angular frequency in the interval $\Tbb:=\left[0,2\pi\right)$, and $\langle \mathbf{k},\boldsymbol{\theta} \rangle :=k_1\theta_1+k_2\theta_2$ is the inner product in $\mathbb{R}^2$. 
Conversely, the covariance $\sigma_\mathbf{k}$ is the $\kb$-th Fourier coefficient of $\Phi(\thetab)$, namely
\begin{equation}\label{moment_eq}
	\sigma_\mathbf{k}=\int_{\mathbb{T}^2}e^{\imath\langle \mathbf{k},\boldsymbol{\theta} \rangle}\Phi(\thetab)\mathrm{d}m(\boldsymbol{\theta}),
\end{equation}
where $\mathrm{d}m(\boldsymbol{\theta})=\frac{1}{{(2\pi)}^2}\mathrm{d}\theta_1\d \theta_2$ is the normalized Lebesgue measure in $\mathbb{T}^2$. 
Therefore, $\Phi$ represents a frequency-domain description of the random field which encodes all the second-order statistics. When the random field is Gaussian, the spectral density in fact specifies the probability distribution of any finite collection of random variables from the field.

The \emph{trigonometric moment equation} \eqref{moment_eq} above plays a fundamental role in a class of spectral estimation methods called \emph{rational covariance extension}. Within that philosophy, one aims to reconstruct a rational spectral density $\Phi$ as a solution to \eqref{moment_eq} given a finite number of covariances $\{\sigma_{\kb}\}_{\kb\in\Lambda}$ where $\Lambda\subset\Zbb^2$ is an index set. 
A typical choice of the index set is simply a ``rectangle'' which will be used throughout this paper:
\begin{equation}\label{set_Lambda}
	\Lambda = \{\kb=(k_1,k_2)\in\Zbb^2 : |k_j| \leq n_j,\ j=1,2\},
\end{equation}
where $n_1$ and $n_2$ are positive integers.
In practice we often observe a finite realization (samples) of the random field
\begin{equation}
	\{ y_{\tb} : t_j\in [0, T_j-1]\cap\Zbb,\ j=1,2\}
\end{equation}
where positive integers $T_1$ and $T_2$ specify the sample size, and the covariances $\{\sigma_{\kb}\}_{\kb\in\Lambda}$ are estimated from the realization with standard formulas, see \citet[p.~23]{stoica2005spectral} and \cite{ZFKZ2019M2}. 
These estimates are reliable when each positive integer $n_j$ is significantly smaller than $T_j$.

Our specific problem setup is taken from \cite{Liu-Zhu-2021} which is a specialization of the work in \cite{ZFKZ2019M2,ZFKZ2020M2-SIAM} by fixing the dimension $d=2$. More precisely, we set up the following \emph{primal} optimization problem
\begin{equation}\label{opt_primal}
	\begin{aligned}
		\min_{\Phi > 0} &\ 
			 D(\Phi,\Psi) := \int_{\mathbb{T}^2} \left[ \log (\Phi^{-1}\Psi) + \Psi^{-1}(\Phi-\Psi) \right] \mathrm{d}m
		\\
		\text{s.t.} & \ \sigma_\mathbf{k}=\int_{\mathbb{T}^2}e^{\imath\langle \mathbf{k},\boldsymbol{\theta} \rangle}\Phi\, \mathrm{d}m
		\quad  \forall \kb \in \Lambda,
	\end{aligned}
\end{equation}
where the Itakura-Saito (IS) pseudo-distance \citep{FMP-12, LP15} is utilized as the objective functional and moment equations are included as equality constraints. The additional element in the problem, $\Psi$, called \emph{prior}, contains some \emph{a priori} information about the desired solution so that one tries to find a $\Phi$ as close as possible to $\Psi$ according to the IS pseudo-distance that meets the moment constraints.
Notice that in \eqref{opt_primal} and the remaining part of this paper, the integral variable $\thetab$ is omitted when there is no risk of confusion.

The primal problem \eqref{opt_primal} is solved via duality. More precisely, the dual optimization problem is written as
\begin{equation}\label{opt_dual}
	\underset{\Qb\in\Lscr_+}{\text{min}} \ J_\Psi(\Qb):=\langle\Qb,\Sigmab\rangle-\int_{\Tbb^2}\log(\Psi^{-1}+Q)\d m
\end{equation}
where, the dual variable $\Qb=\{q_\kb\}_{\kb\in\Lambda}$ consists of the Lagrange multipliers such that $q_{-\kb} = \overline{q_\kb}$ for all $\kb\in\Lambda$ and $q_\zerob\in\Rbb$, $\Sigmab=\{\sigma_\kb\}_{\kb\in\Lambda}$ contains all the covariances in the equality constraints of \eqref{opt_primal}, $\innerprod{\Qb}{\Sigmab}:=\sum_{\kb\in\Lambda} q_\kb \overline{\sigma_{\kb}}$ defines a real-valued inner product between two symmetric multisequences indexed in $\Lambda$, ${Q}(\thetab):=\sum_{\kb \in \Lambda} {q}_{\kb} e^{-\imath\innerprod{\kb}{\thetab}}$ is a symmetric trigonometric polynomial with coefficients in $\Qb$,
and the feasible set is given by
\begin{equation}\label{feasible_set}
	\Lscr_+:=\left\{ \Qb : (\Psi^{-1}+Q)(\thetab)>0 \ \, \forall \thetab\in\Tbb^2 \right\}.
\end{equation}
We remark that $\Qb$ dwells in a vector space of \emph{real} dimension $|\Lambda|=(2n_1+1)(2n_2+1)$ where $|\cdot|$ denotes the cardinality of a set.
The dual problem \eqref{opt_dual} is strictly convex \cite{ZFKZ2019M2}.
Given the optimal solution $\hat{\Qb}$ to the dual problem, the optimal solution to the primal problem can be recovered as $\hat{\Phi}=(\Psi^{-1} + \hat{Q})^{-1}$.



\section{Algebraic structure of the full Hessian}\label{sec:struct_Hessian}



In this section, we investigate the proper algebraic structure of the full Hessian that will enable the true Newton's method to solve the dual problem \eqref{opt_dual}.
Since the specific index set $\Lambda$ in \eqref{set_Lambda} is used in our optimization problem, we can identify the dual variable as a $(2n_2+1)\times(2n_1+1)$ complex matrix
\begin{equation}\label{mat_Q}
	\mathbf{Q} = \bmat
	q_{-n_1,-n_2} & \cdots & q_{0,-n_2} & \cdots & q_{n_1,-n_2} \\
	\vdots & & \vdots & & \vdots \\
	q_{-n_1,0} & \cdots & q_{0,0} & \cdots & q_{n_1,0} \\
	\vdots & & \vdots & & \vdots \\
	q_{-n_1,n_2} & \cdots & q_{0,n_2} & \cdots & q_{n_1,n_2} \\
	\emat.
\end{equation}
%
In order to properly express the gradient and Hessian of the dual function $J_\Psi$, 
let us introduce the usual \emph{vectorization} of the matrix $\Qb$ above by stacking the columns one by one, that is, column $1$ on the top, column $2$ below column $1$, column $3$ below column $2$, etc., and denote the operation as
\begin{equation}
	\qb := \vtr(\Qb).
\end{equation}
The vector $\qb$ has a size of $|\Lambda|$, and $q_{k_1, k_2}$ occupies the $[(2n_2+1)(k_1+n_1)+k_2+n_2+1]$-th entry. 
In fact, this corresponds to the \emph{lexicographic} order of the index set $\Lambda$.
Moreover, we write $[\qb]_j$ for the $j$-th entry of $\qb$.
With an abuse of notation, we shall rewrite the dual function in \eqref{opt_dual} as $J_\Psi(\qb)$ and define the \emph{gradient} as the $|\Lambda|$-dimensional vector 
\begin{equation}
	\nabla J_\Psi(\qb) :=  \frac{\partial J_\Psi(\qb)}{\partial {\qb}}\ \text{such that} \
	\left[\nabla J_\Psi(\qb) \right]_j = \frac{\partial J_\Psi(\qb)}{\partial {[\qb]_j}}.
\end{equation}
For each $j=1, \dots, |\Lambda|$, the last term in the above formula admits the expression
\begin{equation}\label{partial_J_dual}
		\frac{\partial J_\Psi(\qb)}{\partial {q_\kb}} 
		= \sigma_{-\kb} - \int_{\Tbb^2} e^{-\imath\innerprod{\kb}{\thetab}} (\Psi^{-1} +Q)^{-1} \d m
\end{equation}
for a unique $\kb\in\Lambda$, see \cite{ZFKZ2019M2}.
This is known as the \emph{Wirtinger derivative} \citep{hormander1973introduction} of $J_\Psi$
in which the conjugate pair $q_\kb$ and $\overline{q_{\kb}}$ are viewed as formally independent variables in calculations. 

	The second-order partial derivatives can be computed in a similar fashion, and the formula is given by
		\begin{equation}\label{sec_partials_J_dual}
			\begin{aligned}
				\frac{\partial^2 J_\Psi(\qb)}{\partial \overline{q_{\lb}} \partial {q_\kb}} & := \frac{\partial}{\partial \overline{q_{\lb}}} \left[ \frac{\partial J_\Psi(\qb)}{\partial {q_\kb}} \right] \\
				& = \int_{\Tbb^2} e^{\imath\innerprod{\lb-\kb}{\thetab}} (\Psi^{-1}+Q)^{-2} \d m .
			\end{aligned}
		\end{equation}
	It is then natural to organize the second-order partials into a $|\Lambda|\times |\Lambda|$ \emph{Hessian} matrix 
	\begin{equation}\label{Hess_mat}
		\Hcal(\qb) := \frac{\partial}{\partial \overline\qb }\left[\frac{\partial J_\Psi(\qb)}{\partial {\qb}}\right]\ \text{with}\ \left[\Hcal(\qb)\right]_{jk} = \frac{\partial}{\partial \overline{[\qb]_k} }\left[\frac{\partial J_\Psi(\qb)}{\partial {[\qb]_j}}\right]
	\end{equation}
	where each entry $\left[\Hcal(\qb)\right]_{jk}$  correspond to \eqref{sec_partials_J_dual} for a unique pair $(\kb, \lb)\in\Lambda^2$.
	The Hessian $\Hcal(\qb)$ is shown to be (Hermitian) positive definite \cite[Sec.~3.1]{ZFKZ2019M2}.
	The partial Hessian described in \citet[Eq.~(31)]{Liu-Zhu-2021} is in fact a lower-right principal submatrix of $\Hcal(\qb)$ of size $n_2+1+n_1(2n_2+1)\approx|\Lambda|/2$.

We observe that the second-order partials depend on the difference $\lb-\kb$ of the indices. Hence for a feasible $\qb$ and any $\kb\in\Zbb^2$, we can define a complex-valued function
\begin{equation}
	h_\kb(\qb) := \int_{\Tbb^2} e^{\imath\innerprod{\kb}{\thetab}} (\Psi^{-1}+Q)^{-2} \d m.
\end{equation}
Then the Hessian can be viewed as a  Hermitian block-Toeplitz matrix
\begin{subequations}\label{H-matrix}
	\begin{equation}\label{H-Q}
		\Hcal  = \bmat 
		\mathbf{R}_0 & \mathbf{R}_1 & \dots & \mathbf{R}_{2n_1}\\
		\mathbf{R}_{1}^\ast & \ddots & \ddots & \vdots\\
		\vdots & \ddots & \ddots & \mathbf{R}_1\\
		\mathbf{R}_{2n_1}^\ast & \dots & \mathbf{R}_{1}^\ast & \mathbf{R}_0 
		\emat 
	\end{equation}
	with $2n_1+1$ blocks in each row (and column)
	 where each block
	\begin{equation}\label{R-j}
		\mathbf{R}_j = \bmat 
		h_{j,0} & h_{j,1} & \dots & h_{j,2n_2}\\
		h_{j,-1} & \ddots & \ddots & \vdots\\
		\vdots & \ddots & \ddots & h_{j,1}\\
		h_{j,-2n_2} & \dots & h_{j,-1} & h_{j,0} 
		\emat
	\end{equation}
\end{subequations}
is a square  Toeplitz matrix of size $2n_2+1$.
For simplicity, we have omitted the dependence on $\qb$ in \eqref{H-matrix}. Such a structure of $\Hcal$ is called Toeplitz-block Toeplitz in \cite{wax1983efficient}, and TBT for short.
It follows that in practical calculations, only the first  block row of the Hessian needs to be computed.


\section{Fast inversion of the structured Hessian}\label{sec:fast_inv}

To simplify the discussion, let us consider the case $n_1=n_2=n$ for some positive integer $n$, and take $p=2n+1$ for the size of each block in \eqref{R-j}. We can now define a sequence of matrices $\{\Hcal_j\}_{j\geq 0}$ of increasing size such that $\Hcal_0=\Rb_0$, and 
\begin{equation}\label{H_j}
	\Hcal_j  = \bmat 
	\mathbf{R}_0 & \mathbf{R}_1 & \dots & \mathbf{R}_{2j}\\
	\mathbf{R}_{1}^\ast & \ddots & \ddots & \vdots\\
	\vdots & \ddots & \ddots & \mathbf{R}_1\\
	\mathbf{R}_{2j}^\ast & \dots & \mathbf{R}_{1}^\ast & \mathbf{R}_0 
	\emat \quad\text{for}\ j\geq 1.
\end{equation}
Clearly $\{\Hcal_j\}_{j\geq 0}$ has a nested structure in the sense that $\Hcal_j$ is a principal submatrix of $\Hcal_{j+1}$ for all $j\geq 0$.
A more precisely statement can be made as follows.
Let us define
\begin{equation}\label{underline_R}
	\begin{aligned}
		\underline{\mathbf{R}}_{k} & := \bmat
		\mathbf{R}_{k} \\
		\vdots \\
		\mathbf{R}_1
		\emat,\quad k=1, 2, \dots,
	\end{aligned}
\end{equation}
and an accompanying sequence of matrices of increasing size $\{\Gcal_k\}_{k\geq 0}$ by  $\Gcal_0=\Rb_0$ and
\begin{equation}\label{H_n+1_1}
	\begin{aligned}
		\Gcal_{k} = \bmat
		\Gcal_{k-1} & \underline{\mathbf{R}}_{k} \\
		\underline{\mathbf{R}}_{k}^\ast & \mathbf{R}_0 
		\emat,\quad k=1, 2, \dots.
	\end{aligned}
\end{equation}
Then we have
\begin{equation}\label{H_n+1_2}
		\Hcal_{j} = \Gcal_{2j},\quad j=1, 2, \dots.
\end{equation}
The nested structure is reflected in \eqref{H_n+1_1}.

Fast inversion of TBT matrices has been studied in \cite{wax1983efficient}, and the main idea is to invert a sequence of nested TBT matrices $\{\Gcal_k\}_{k\geq 0}$, using the computed inverse of $\Gcal_{k-1}$ to facilitate the computation of $\Gcal_k^{-1}$.
The algorithm is briefly reviewed next.
Notice that our $\Hcal_{j}$ and $\Gcal_k$ are Hermitian, which can bring about additional saving.
Also, the inversion of our Hessian \eqref{H-Q}, which is now $\Hcal_{n}$ in view of \eqref{H_j}, follows an identical procedure because of \eqref{H_n+1_2}. 


A fundamental property of TBT matrices is that they are symmetric with respect to the antidiagonal, i.e., the diagonal connecting the lower-left and upper-right corners of a square matrix. Such a property is called \emph{persymmetry} in \citet[Lemma~1]{wax1983efficient}, and is described mathematically by the formula
\begin{equation}\label{TBT_persym}
	\Gcal_k = \mathit{J}_{p,k+1} \Gcal_k^{\top} \mathit{J}_{p,k+1},
\end{equation}
where  $\mathit{J}_{p,k}$ is the $pk\times pk$ permutation matrix \citep{zohar1969toeplitz}
\begin{subequations}
	\begin{equation}\label{J_pk}
		\begin{aligned}
			\mathit{J}_{p,k} = \bmat 
			&                                        & \mathit{J}_p \\
			& \begin{sideways}$\ddots$\end{sideways} & \\
			\mathit{J}_p &                                        & 
			\emat
		\end{aligned}
	\end{equation}
	with $k\times k$ blocks
	and each block $\mathit{J}_p$ is the $p\times p$ permutation matrix
	\begin{equation}
		\begin{aligned}
			\mathit{J}_p = \bmat 
			&                                        & 1 \\
			& \begin{sideways}$\ddots$\end{sideways} & \\
			1 &                                        &
			\emat .
		\end{aligned}
	\end{equation}
\end{subequations}

Using block elimination, the block matrix in \eqref{H_n+1_1} admits a block inverse
\begin{subequations}\label{R1}
\begin{equation}\label{G_k_inv_1}
	\begin{aligned}
		\Gcal_{k}^{-1} = \bmat
		\Gcal_{k-1}^{-1} + \underline{\mathit{W}}_{k}\boldsymbol{\alpha}_k^{-1}\underline{\mathit{W}}_{k}^\ast & \underline{\mathit{W}}_{k}\boldsymbol{\alpha}_k^{-1} \\
		\boldsymbol{\alpha}_k^{-1}\underline{\mathit{W}}_{k}^\ast & \boldsymbol{\alpha}_k^{-1}
		\emat
	\end{aligned}
\end{equation}	
where
\begin{equation}\label{W_underline_k}
	\begin{aligned}
		\underline{\mathit{W}}_{k} := -\Gcal_{k-1}^{-1}\underline{\mathbf{R}}_{k} 
	\end{aligned}
\end{equation}
and
\begin{equation}\label{alpha_k}
	\begin{aligned}
		\boldsymbol{\alpha}_k := \Rb_0 + \underline{W}_k^* \underline{\Rb}_{k}
	\end{aligned}
\end{equation}
is the \emph{Schur complement $\Gcal_{k}/\Gcal_{k-1}$}.
\end{subequations}
Now let us take inverse of both sides of the equation in \eqref{TBT_persym}, and we have
\begin{equation}\label{R2}
\begin{aligned}		
	\Gcal_{k}^{-1} & = 
	\mathit{J}_{p,k+1} \Gcal_k^{-\top} \mathit{J}_{p,k+1} \\
	& = \bmat
	\mathit{J}_p \boldsymbol{\alpha}_{k}^{-\top} \mathit{J}_p & \mathit{J}_p \boldsymbol{\alpha}_{k}^{-\top} \underline{\widehat{\mathit{W}}}_{k}^\top \\
	\overline{\underline{\widehat{\mathit{W}}}}_{k} \boldsymbol{\alpha}_{k}^{-\top} \mathit{J}_p & \Gcal_{k-1}^{-1} + \overline{\underline{\widehat{\mathit{W}}}}_{k} \boldsymbol{\alpha}_{k}^{-\top} \underline{\widehat{\mathit{W}}}_{k}^\top 
	\emat ,
\end{aligned}
\end{equation}
where, the second equality follows from direct computation with \eqref{J_pk} and \eqref{G_k_inv_1}, $\overline{A}$ denotes elementwise complex conjugate of a matrix $A$, and
\begin{equation}\label{def_hat_W_underline}
	\begin{aligned}
		\underline{\widehat{W}}_{k} := \mathit{J}_{p,k}\underline{\mathit{W}}_k 
	\end{aligned}
\end{equation}
with $\underline{W}_k$ in \eqref{W_underline_k}.

From the two representation of $\Gcal_k^{-1}$  \eqref{G_k_inv_1} and \eqref{R2}, we can read out two equations between the blocks
\begin{equation}\label{H_ij}
\begin{aligned}
	\left( \Gcal_{k}^{-1} \right)_{i, j} = \left( \Gcal_{k-1}^{-1} + \underline{W}_{k}\boldsymbol{\alpha}_k^{-1}\underline{W}_{k}^\ast\right)_{i, j},\quad i,j = 1,\dots, k,
\end{aligned}
\end{equation}
and
\begin{equation}\label{H_ij+1}
	\left( \Gcal_{k}^{-1} \right)_{i+1, j+1} = \left( \Gcal_{k-1}^{-1} + \overline{\underline{\widehat{W}}}_{k}\boldsymbol{\alpha}_k^{-\top}\underline{\widehat{W}}_{k}^\top \right)_{i, j}, \quad i,j = 1,\dots, k,
	\end{equation} 
	where the subscripts $i$ and $j$ refer to blocks.
	Combining \eqref{H_ij} and \eqref{H_ij+1}, the blocks of $\Gcal_{k-1}^{-1}$ can be eliminated, and we are left with a recursion among blocks of $\Gcal_k^{-1}$:
	\begin{equation}\label{recursion_block_G_k_inv}
				\left(\Gcal_{k}^{-1} \right)_{i+1, j+1} = \left(\Gcal_{k}^{-1} \right)_{i, j} +
				(\Wcal_k)_{i, j},\quad i, j=1, \dots, k. 
	\end{equation}
	where $\Wcal_k:=\overline{\underline{\widehat{W}}}_{k}\boldsymbol{\alpha}_k^{-\top}\underline{\widehat{W}}_k^\top - \underline{W}_{k}\boldsymbol{\alpha}_k^{-1}\underline{W}_{k}^\ast$ is a $k\times k$ block matrix.
	If we introduce the new notation
	\begin{subequations}
		\begin{equation}\label{def_hat_R}
			\widehat{\mathbf{R}}_j := \mathit{J}_p \mathbf{R}_j,\quad j=1, 2, \dots,
		\end{equation}
		and
		\begin{equation}\label{def_hat_R_underline}
			\underline{\widehat{\mathbf{R}}}_k := \mathit{J}_{p, k} \underline{\mathbf{R}}_k = 
			\bmat
			\widehat{\mathbf{R}}_1 \\ \vdots \\\widehat{\mathbf{R}}_k
			\emat ,\quad k=1, 2, \dots,
		\end{equation}
	\end{subequations}
	where $\underline{\Rb}_k$ was defined in \eqref{underline_R},
%
then according to Formulas (13a), (13c), and (17) in \cite{wax1983efficient}, the recursions for $\underline{\widehat{W}}_{k}$ and $\boldsymbol{\alpha}_k$ are given as follows:
\begin{subequations}\label{W,V,alpha}
\begin{equation}\label{W}
	\begin{aligned}
		\underline{\widehat{\mathit{W}}}_{k+1} = \bmat
		\underline{\widehat{W}}_{k} \\
		\mathbf{0}_{p}
		\emat
		- \bmat
		\underline{\overline{\mathit{W}}}_{k} \\
		\mathit{I}_{p}
		\emat 
		\boldsymbol{\alpha}_k^{-\top}\boldsymbol{\beta}_k, 
	\end{aligned}
\end{equation}
and
\begin{equation}\label{alpha}
	\begin{aligned}
		\boldsymbol{\alpha}_{k+1} = \boldsymbol{\alpha}_k - \boldsymbol{\beta}_k^\ast \boldsymbol{\alpha}_k^{-\top} \boldsymbol{\beta}_k,
	\end{aligned}
\end{equation}
where
\begin{equation}\label{beta,gamma,mu}
	\begin{aligned}
		\boldsymbol{\beta}_k = \underline{W}_k^\top \underline{\widehat{\mathbf{R}}}_{k} + \widehat{\mathbf{R}}_{k+1} . 
	\end{aligned}
\end{equation}
\end{subequations}

The pseudocode for the inversion of our structured Hessian is given in Algorithm~\ref{alg:Inversion}.
Notice that due to the Hermitian and persymmetric structure of $\Gcal_k^{-1}$, roughly only a quarter of the blocks need to be computed, and this point is also illustrated in Algorithm~\ref{alg:Inversion}.
\begin{algorithm}
\begin{algorithmic}[1]
	\caption{Inversion of the Hessian $\Hcal_n=\Gcal_{2n}$}\label{alg:Inversion}
	\Require A Hessian matrix of the form \eqref{H-matrix}.
	\State $\Gcal_0^{-1} \gets \mathbf{R}_0^{-1}$.
	\State Compute $\underline{\mathit{W}}_1$,  $\boldsymbol{\alpha}_1$ 
    by \eqref{R1}.
	\State Compute $\underline{\widehat{\mathit{W}}}_1$, $\underline{\widehat{\mathbf{R}}}_1$ by \eqref{def_hat_W_underline} and \eqref{def_hat_R_underline}, respectively.
	\While{$1\leq k\leq 2n-1$}
	\State Compute $\widehat{\mathbf{R}}_{k+1}$ by \eqref{def_hat_R}.
	\State Update $\boldsymbol{\beta}_k$ by \eqref{beta,gamma,mu}.
	\State Update $\boldsymbol{\alpha}_{k+1}$ by \eqref{alpha}.
	\State Update $\underline{\widehat{\mathit{W}}}_{k+1}$ by \eqref{W}.
	\EndWhile
	\State $\Gcal_{2n}^{-1}\gets p^2 \times p^2$ zero matrix.
	\State Update the first block row of $\Gcal_{2n}^{-1}$ by 
	\begin{equation*}
			\left(\Gcal_{k}^{-1} \right)_{1, 1}  = J_p \boldsymbol{\alpha}_{k}^{-\top} J_p \ \ \text{and} \ \
			\left(\Gcal_{k}^{-1} \right)_{1, 2:k+1} = \mathit{J}_p \boldsymbol{\alpha}_{k}^{-\top} \underline{\widehat{W}}_{k}^\top.
	\end{equation*}
	\State Recover $\underline{\mathit{W}}_{2n} = \mathit{J}_{p,2n} \underline{\widehat{W}}_{2n}$.
	\State // \texttt{Update the upper quarter of $\Gcal_{2n}^{-1}$.}
	\For{$i=1:n$}
	\For{$j=i:2n-i$}
	\State Update $\left( \Gcal_{2n}^{-1} \right)_{i+1, j+1}$ by \eqref{recursion_block_G_k_inv}.
	\EndFor 
	\EndFor 
	\State // \texttt{Halve the diagonal and antidiagonal elements of $\Gcal_{2n}^{-1}$. Notice that the center block $\left(\Gcal_{2n}^{-1}\right)_{n+1, n+1}$ gets quartered.}
	\For{$i=1:n+1$}
	\State $\left( \Gcal_{2n}^{-1} \right)_{i, i} \gets \frac{1}{2}\left( \Gcal_{2n}^{-1} \right)_{i, i}$.
	\State $\left( \Gcal_{2n}^{-1} \right)_{i, 2n+2-i} \gets \frac{1}{2}\left( \Gcal_{2n}^{-1} \right)_{i, 2n+2-i}$.
	\EndFor 
	\State // \texttt{Update the right quarter of $\Gcal_{2n}^{-1}$ by persymmetry.}
	\State $\Gcal_{2n}^{-1}\gets \Gcal_{2n}^{-1} + \left(\mathit{J}_{p,2n+1} \Gcal_{2n}^{-1} \mathit{J}_{p,2n+1}\right)^\top$.
	\State // \texttt{Update the block lower triangular part  of $\Gcal_{2n}^{-1}$ by Hermitianity.}
	\State $\Gcal_{2n}^{-1}\gets \Gcal_{2n}^{-1} + \left(\Gcal_{2n}^{-1}\right)^*$.
	\\
	\Return $\Gcal_{2n}^{-1}$.
\end{algorithmic}
\end{algorithm}


\section{Fast solution to TBT systems of linear equations}\label{sec:TBT_linear_eq}

Fast solution to TBT systems of linear equations was also studied in \cite{wax1983efficient}. Here we
slightly expand their results to the case in which the unknown and the right-hand matrices can have an arbitrary number of columns. To this end, consider the linear equations
\begin{equation}\label{linear equations}
	\begin{aligned}
		\Hcal_{n} \underline{X}_{2n}=\underline{B}_{2n}
	\end{aligned}
\end{equation}
where $\Hcal_{n}$ is the TBT Hessian matrix in \eqref{H-Q}, and $\underline{X}_{2n}$ (unknown) and $\underline{B}_{2n}$ (known) are $p^2\times q$ matrices. Apparently, $\underline{X}_{2n}$ and $\underline{B}_{2n}$ can be viewed as $p\times 1$ block vectors where each block has a size of $p\times q$.

Note that due to the relation \eqref{H_n+1_2} and the nested structure \eqref{H_n+1_1}, we can define a sequence of equations
\begin{equation}\label{linear question}
	\begin{aligned}
		\Gcal_{k}\underline{X}_k=\underline{B}_k, \quad   0\leq k\leq 2n,
	\end{aligned}
\end{equation}
where $\Gcal_{k}$ has been defined in \eqref{H_n+1_1}, $\underline{B}_k$ is formed by the leading $k+1$ blocks of the block vector $\underline{B}_{2n}$, and $\underline{X}_k$ is an unknown matrix of size $(k+1)p\times q$. Notice that in general $\underline{X}_k$ is not equal to the leading $k+1$ blocks of $\underline{X}_{2n}$, and only when $k=2n$ we obtain the solution $\underline{X}_{2n}$ to \eqref{linear equations}.

Using \eqref{R1} we can derive the formula
\begin{subequations}\label{linear answer}
	\begin{equation}\label{X_k+1}
		\begin{aligned}
			\underline{X}_{k+1} =  \bmat
			\underline{X}_k \\
			\mathbf{0}_{p\times q}
			\emat
			+ \bmat
			\underline{{\mathit{W}}}_{k} \\
			\mathit{I}_p
			\emat 
			\boldsymbol{\alpha}_{k}^{-1}\boldsymbol{\delta}_{k}
		\end{aligned}
	\end{equation}
	and
	\begin{equation}\label{delta_k}
		\begin{aligned}
			\boldsymbol{\delta}_k = \underline{\mathit{W}}_k^\ast \underline{B}_k+B_{k+1}, \quad   0\leq k\leq 2n-1,
		\end{aligned}
	\end{equation}
\end{subequations}
where $B_{k}$ (without underline) is the $(k+1)$-th block of the block vector $\underline{B}_{2n}$, and $\underline{{\mathit{W}}}_{k}$ and $\boldsymbol{\alpha}_{k}$ are updated via \eqref{def_hat_W_underline}, \eqref{W}, and \eqref{alpha}.
We summarize the pseudocode for the fast solution to TBT systems of linear equations in Algorithm~\ref{alg:linearsolver}.

\begin{algorithm}
	\begin{algorithmic}[1]
		\caption{Solving the TBT system \eqref{linear equations} of equations}\label{alg:linearsolver}
		\Require A Hessian matrix of the form \eqref{H-Q} and a $p^2\times q$ matrix $\underline{B}_{2n}$.
		\State $\bmat\underline{X}_0 & \underline{W}_1\emat \gets \mathbf{R}_0\setminus \bmat\underline{B}_0 & -\Rb_1\emat$.
		\State Compute 
        $\boldsymbol{\alpha}_1$ by \eqref{alpha_k}.
		\State Compute $\underline{\widehat{\mathit{W}}}_1$, $\underline{\widehat{\mathbf{R}}}_1$ by \eqref{def_hat_W_underline} and \eqref{def_hat_R_underline}, respectively.
		\While{$1\leq k\leq 2n-1$}
		\State Update $\boldsymbol{\delta}_k$ by \eqref{delta_k}.
		\State Update $\underline{X}_k$ by \eqref{X_k+1}.
		\State Compute $\widehat{\mathbf{R}}_{k+1}$ by \eqref{def_hat_R}.
		\State Update $\boldsymbol{\beta}_k$ by \eqref{beta,gamma,mu}.
		\State Update $\boldsymbol{\alpha}_{k+1}$ by \eqref{alpha}.
		\State Update $\underline{\widehat{\mathit{W}}}_{k+1}$ by \eqref{W}.
		\EndWhile
		\State Compute $\boldsymbol{\delta}_{2n}$ by \eqref{delta_k}.
		\State Compute $\underline{X}_{2n}$ by \eqref{X_k+1}.
		\\
		\Return $\underline{X}_{2n}$.
	\end{algorithmic}
\end{algorithm}

\section{Simulations}\label{sec:sims}

In this section, we perform numerical simulations for the fast inversion of positive definite TBT matrices and for the solution of the dual problem \eqref{opt_dual} using Newton's method with the full Hessian.
The results are divided into two subsections.

\subsection{Testing Algorithm~\ref{alg:Inversion}}

\begin{figure}[htb]
	\begin{center}
		\includegraphics[width=6.4cm]{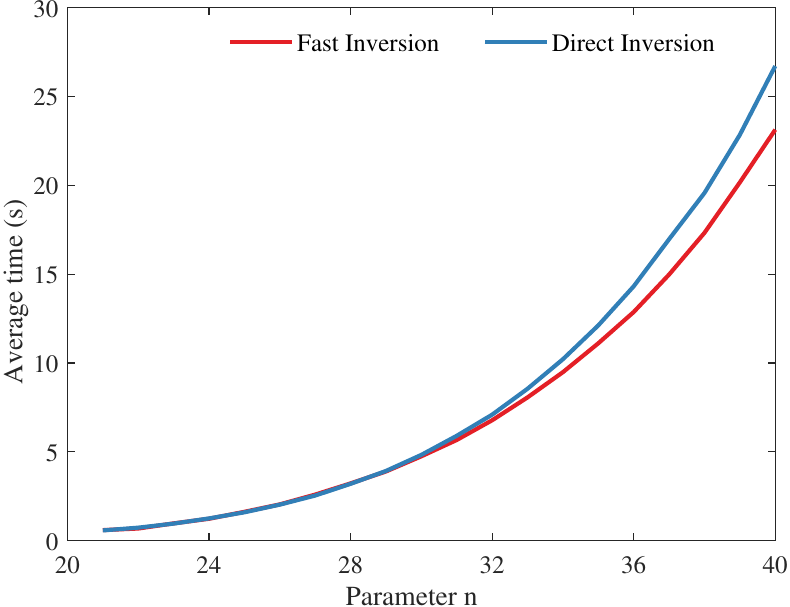}    
		\caption{Algorithm~\ref{alg:Inversion} vs. the direction method for the inversion of positive definite TBT matrices: average running times for $n$ ranging from $21$ to $40$.}
		\label{fig:Average-time}
	\end{center}
\end{figure}

The first set of numerical experiments are about  the inversion of positive definite TBT matrices. In accordance with the previous sections, the size parameters  are set as $n_1 = n_2 = n$, resulting in $p^2\times p^2$ TBT matrices of the form \eqref{H-matrix} with $p=2n+1$.
Each TBT matrix is randomly generated and is made positive definite by adding a positive multiple of the identity matrix. Then
two methods are compared for the inversion of such a TBT matrix, namely Algorithm~\ref{alg:Inversion} for fast inversion and the direct inversion using the MATLAB command \verb|A^(-1)| for any square matrix $\texttt{A}$. We compute these inverses for values of $n$ ranging from $1$ to $40$. Each each $n$, twenty positive definite TBT matrices are generated and the running times of the two methods are averaged and shown in Fig.~\ref{fig:Average-time}.
Notice that we only present the running times for $20< n\leq 40$, because for $n\leq 20$, both methods yield nearly identical results.
We see that Algorithm~\ref{alg:Inversion} outperforms the direction inversion when matrix size is large, in particular, when $n\geq 30$. 
The advantage of the fast algorithm can be explained theoretically.
Indeed, according to \cite{wax1983efficient}, the computational complexity of Algorithm~\ref{alg:Inversion} is $O(p^5)$, see also \cite{Akaike-1973-BT-inv}, while the direct inversion has a complexity of $O(p^6)$, the usual cubic complexity applied to a square matrix of size $p^2$.



\subsection{Full Hessian vs. quarter Hessian in Newton's method}

In this subsection, we use Newton's method with the full Hessian to solve the dual problem \eqref{opt_dual}, and make a comparison with the quasi-Newton method in \cite{Liu-Zhu-2021} which involves only a quarter Hessian, see the discussion after \eqref{Hess_mat}.


In our implementation, all integrals against the normalized Lebesgue measure such as those in \eqref{opt_dual}, \eqref{partial_J_dual}, and \eqref{sec_partials_J_dual} are discretized as Riemann sums:
\begin{equation}
	\int_{\mathbb{T}^2}f \mathrm{d}m \approx \frac{1}{N_1N_2}\sum_{\thetab \in \mathbb{T}_\mathbf{N}^2}f(\thetab),
\end{equation}
where
	$\mathbb{T}_\mathbf{N}^2:=\left\{\left( \frac{2\pi}{N_1}\ell_1,\frac{2\pi}{N_2}\ell_2\right) : \ell_j \in [0, N_j-1],\ j=1,2 \right\}$
is a regular grid for $\Tbb^2$ defined by positive integers $N_1$ and $N_2$.
The 2-D random field under consideration admits a model of a single complex exponential corrupted by additive white complex Gaussian noise, see Sec.~5 in \cite{Liu-Zhu-2021}. The ratio between the amplitude of the complex exponential and the standard deviation of the noise is $1/\sqrt{2}$. The sample size of the random field is $N_1\times N_2$, the same as the size of the grid $\Tbb^2_{\Nb}$. Standard biased covariance estimates $\{\hat \sigma_\kb\}_{\kb\in\Lambda}$ \citep{stoica2005spectral} are computed from the samples and they play the role of
$\Sigmab=\{\sigma_\kb\}_{\kb\in\Lambda}$ in the dual problem. In addition, the prior $\Psi\equiv \sigma_{\zerob}$ is a constant.



\begin{figure}[htb]
	\begin{center}
		\includegraphics[width=6.4cm]{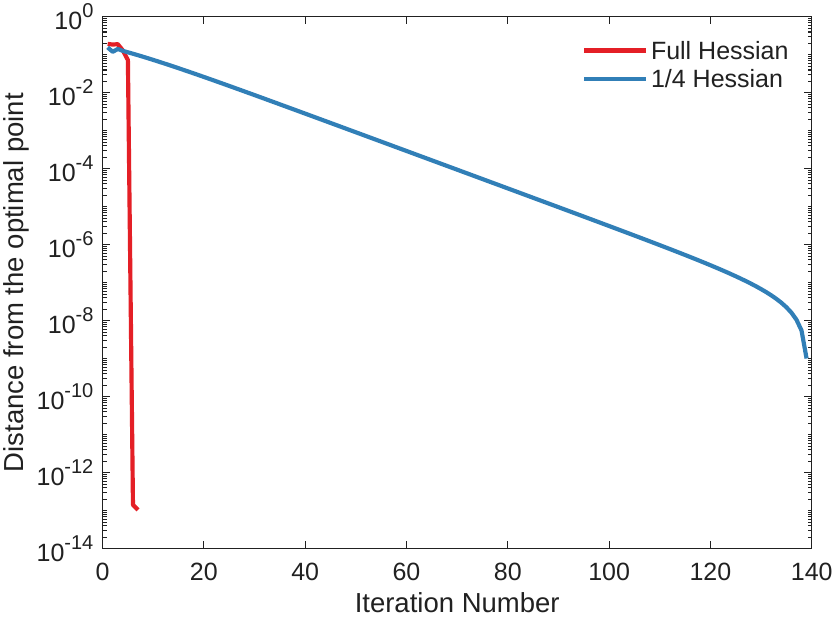}    
		\caption{Full Hessian vs. the quarter Hessian in Newton's method for the solution of the dual problem \eqref{opt_dual}: Euclidean distance between the $k$-th iterate $\qb^{k}$ and the optimal point $\qb^*$. }
		\label{fig:distance}
	\end{center}
\end{figure}

In the numerical example to be presented next, we take $N_1=N_2=30$, and $n_1=n_2=1$ which defines the index set $\Lambda$ in \eqref{set_Lambda}. 
Newton's method involves the full  Hessian of size $9\times9$, and the Newton direction at each iteration is obtained by solving a linear system of equations of the form \eqref{linear equations} using Algorithm~\ref{alg:linearsolver}. In contrast, the quasi-Newton method in \cite{Liu-Zhu-2021} employs a reduced Hessian of size $5\times5$.
The convergence of the two methods with respect to the number of iterations are depicted in Fig.~\ref{fig:distance}.
It is clear that the true Newton's method converges much more rapidly than the alternative. 
Although each iteration of the full-Hessian computation takes approximately twice of the time for the quarter Hessian,  overall the computational times using the full Hessian and the quarter Hessian are approximately $0.04$ and $0.38$ seconds, respectively, when the accuracy of the iterate as measured by the distance to the optimal point is within $10^{-8}$.
Thus, in this example Newton's method is faster then the quasi-Newton method using the quarter Hessian by roughly an order of magnitude.

\section{Conclusions}\label{sec:conclus}

We have devised a fast Newton solver for the solution of a 2-D spectral estimation problem. The Toeplitz-block Toeplitz structure of the full Hessian of the dual function is exploited so that an efficient inversion algorithm from the literature can be used upon suitable adaptation. The convergence of the true Newton's method is remarkable in comparison with the quasi-Newton method in a previous work which uses only a quarter Hessian.



\bibliography{ifacconf}
  

\end{document}